\documentclass[12pt]{article} 
\usepackage[utf8]{inputenc}
\usepackage{mathtools}
\usepackage{amstext}
\usepackage{amsfonts}
\usepackage{amssymb}
\usepackage{amsthm}
\usepackage{amsmath}
\usepackage{url}
\usepackage{enumitem}
\usepackage{cases}
\usepackage{verbatim}
\usepackage{xcolor}
\usepackage{footmisc}

\usepackage{mathrsfs, url}
\usepackage[colorlinks=true,linkcolor=blue]{hyperref}
\usepackage[margin=1in]{geometry}

\newcommand{\pr}[1][]{\mathbb{P}}

\def\N{{\mathbb N}}

\newtheorem{thm}{Theorem}[section]

\newtheorem{cor}[thm]{Corollary}
\newtheorem{lemma}[thm]{Lemma}
\newtheorem{prop}[thm]{Proposition}

\newtheorem{defn}[thm]{Definition}

\newtheorem{example}[thm]{Example}

\newcommand{\remove}[1] {}

\DeclareMathOperator{\vol}{vol}
\DeclareMathOperator{\dist}{dist}
\DeclareMathOperator{\diam}{diam}
\DeclareMathOperator{\ecc}{ecc}
\DeclareMathOperator{\dom}{dom}
\DeclareMathOperator{\Dom}{Dom}

\author{Erin Meger \thanks{School of Computing, Queen's University, Kingston, ON (erin.meger@queensu.ca)}
\and
Abigail Raz\thanks{Department of Mathematics, The Cooper Union for the
Advancement of Science and Art, NYC, NY
  (abigail.raz@cooper.edu) - Corresponding author}}

\title{The Iterative Independent Model}

\begin{document}
\maketitle

\begin{abstract}
    Deterministic complex networks that use iterative generation algorithms have been found to more closely mirror properties found in real world networks than the traditional uniform random graph models. 
    In this paper we introduce a new, Iterative Independent Model (IIM), generalizing previously defined models in \cite{ilm, ilt, ilat}. These models use ideas from Structural Balance Theory to generate edges through a notion of cloning where ``the friend of my friend is my friend'' and anticloning where ``the enemy of my enemy is my friend'' \cite{ilm, easley2010networks}.
    In this paper, we vastly generalize these notions by allowing each vertex added at a given time step to choose independently of the other vertices if it will be cloned or anticloned. While it may seem natural to focus on a randomized model, where we randomly determine whether or not to clone any given vertex, we found the general deterministic model exhibited certain structural properties regardless of the probabilities. This allows applications to then explore the particulars, while having the theoretical model explain the structural phenomenons that occur in all possible scenarios. 
    
    Throughout the paper we demonstrate that all IIM graphs have spectral gap bounded away from zero, which indicates the clustering properties also found in social networks. Furthermore, we show bounds on the diameter, domination number, and clique number further indicating the well clustered behaviour of IIM graphs. Finally, for any fixed graph $F$ all IIM graphs will eventually contain an induced copy of $F$.

\end{abstract}

\section{Introduction} 


Models of complex networks are used to predict the evolution of structure within real-world networks, and rely on observed properties and phenomenon to influence generation algorithms including internet traffic, biological interactions, and social networks \cite{aclpa,barabasi2012luck, barabasiPA, chung2003duplication, erdosreyni, estrada2005complex}. For a complete introduction to complex networks, their definitions, and their properties, see \cite{fanbook} and \cite{easley2010networks}.
Recently, there has been a development in deterministic complex networks that use iterative generation algorithms \cite{ilm, ilt, ilat, igm}. These models, in particular, use underlying social network principles to define local interactions reflecting the transitivity of in-group and out-group behaviour \cite{easley2010networks, zachary1977information}. In particular, the Iterated Local Model is generated by creating a new node for each node in the previous iteration, $G$. At each time step, we either clone every node $v \in V(G)$, by adding $v'$, a vertex with precisely $v$'s neighborhood, including $v$ itself, or anticlone every node $v \in V(G)$ adding $v^\star$ whose neighborhood is precisely the vertices not adjacent to $v$ nor $v$ itself. In this way, all nodes experience an in-group expansion ``the friend of my friend is my friend'' or an out-group expansion ``the enemy of my enemy is my friend'' \cite{ilm, easley2010networks}. In this paper we seek to generalize this notion by allowing vertices to experience cloning/anticloning independently.

We define the Independent Iterated Model (IIM) where each node experiences either a transitive cloning or an anti-transitive anticloning, independently of all other nodes added in the same time step. This allows us to study the significantly more general structure that arises from the cloning/anticloning process. While a probabilistic approach seems intuitive, the model exhibits inherent structural properties regardless of any randomness. Thus, these properties are intrinsic within the transitive nature of social networks, and can be more easily applied and studied. That is, for any application, regardless of the specific probability used to generate such an IIM graph, one can immediately apply the results in this paper to determine specific structural properties. Furthermore, the techniques used in this paper demonstrate that these properties rely primarily on the transitive/anti-transitive nature rather than just the simple fact that graphs generated using such a model have large independent sets. 
The complete generalization of this class of complex network models opens the door for a wide array of applications to social network influence, internet traffic, and information diffusion \cite{stoica2020seeding, palu, kempe2003maximizing}. 

This paper is organized as follows. In Section \ref{prelim} we recall some basic graph notation and formally define our model. In Section \ref{spectral} we show that the spectral gap of any graph generated using our model is bounded away from zero. In Section \ref{graphprop} we explore a number of graph properties of our model including diameter, domination number, clique number, and coloring. Finally, in Section \ref{induced} we show that for any fixed graphs $F$ and $G$ there is a sufficiently large $k$ such that our model after $k$ steps originating with $G$ contains $F$ as an induced subgraph.

\section{Preliminaries}\label{prelim}
We begin with a discussion of notation and terminology we will use throughout the paper. Although many notations are common, we recall them here to avoid ambiguity. Then, we define the deterministic generation algorithm for the Iterated Independent Model.
\subsection{Graph Notation}

Given a graph $G$ let $V(G)$ denote the vertex set of $G$. We use $N(v)$ to denote the neighborhood of $v \in V(G)$ and use the terms \emph{closed neighborhood of $v$} for $N[v]= \{v\} \cup N(v)$ and \emph{anti-neighborhood of $v$} for $V(G) \setminus N[v]$. For any set $X \subseteq V(G)$ and $v \in V(G)$ we denote $\deg_X(v)$ for the number of neighbors of $v$ in $X$. For any $X, Y \subseteq V(G)$ let $E(X,Y)$ denote the set of edges of $G$ with one endpoint in $X$ and one end point in $Y$; for the case of a single edge with endpoints $x$ and $y$ we simply write $xy$ for the edge. Additionally, we let $E(X)=E(X,X)$, $n=|V(G)|$, $\overline{X}= V(G) \setminus X$, and say $\vol(X)=  \sum_{v \in X} \deg(v)$.  
For any graph $H$ and $X \subseteq V(X)$ we let $H[X]$ denote the induced subgraph of $H$ on $X$, that is the graph with vertex set $X$ such that for $x,y \in X$ we have the edge $xy$ if and only if $xy \in E(H)$. There are a number of additional graph terms and notation that we will define at the beginning of each relevant section.


\subsection{The IIM Model}\label{model}

As mentioned above in social networks, when presented with many options, we often make choices as to which connections we would like to have. This was a key competent to the analysis of the Zachary Karate club, where club members split along the connections with respect to different instructors \cite{zachary1977information}. In this way, we can consider a model where each new member may decide they liked or disliked their initial instructor independently. Thus we introduce the \emph{Iterated Independent Model (IIM)}.

The Iterated Independent Model is a complex network generated by a deterministic process that grows the graph at each time step following a particular generation algorithm as defined below. Many special cases of our model have been previously studied, the most general of which being the Iterated Local Model \cite{ilm}.

Given an initial graph $G$ and vertex $v \in V(G)$ we say a new vertex $v'$ is $v$'s clone if 
\[N(v')=N[v]
\]
and $v$'s anticlone if
\[N(v')=V(G)\setminus N[v].
\]
In the previous work, the ILM graphs were initiated with some starting graph $G_0$, and at each time step $t \geq 1$ built $G_t$ from $G_{t-1}$ by either cloning every vertex in $G_{t-1}$ or anticloning every vertex in $G_{t-1}$ \cite{ilm}. The ILM itself generalizes two previous models, the Iterated Local Transitive Model and the Iterated Local Anti-Transitive Model, which each used exclusively the cloning or anticloning procedures, respectively \cite{ilt, ilat}. These models exhibit similar structural properties, which we show are fully generalizable via our Iterated Independent Model.

The IIM initiates with some graph $G_0$. At each time step we independently choose, for each previously existing vertex, whether to clone or anticlone as in the definitions above. Each vertex behaves independently of other vertices, and the structures that more easily appear in ILM are no longer necessarily present. Interestingly, we found many similar properties still hold but new techniques are needed. Many results from \cite{ilm} hold as special cases of the more generalized results presented in this paper. This allows for generalizations to further applications, to consider many more complex scenarios.

We use the notation $\mathcal{IIM}_l(G)$ to indicate the set of all possible graphs generated using this model originating with graph $G$ and allowed to run for $l$ time steps. Note we say $\{G\}=IIM_0(G)$. 
We will refer to the vertices added in each time step as a \emph{level}, and note that $V(G)$ is level 0. If $v$ is either a clone or an anticlone of $u$ then we say $v$ is a \emph{copy} of $u$ and $u$ is a \emph{precopy} of $v$. If we have some set $\{v_0, \ldots, v_k\} \subseteq V(G_l)$ such that each $v_{i+1}$ is a copy of $v_i$, for $i \ge 0$ then we say that for any $i \ge j$ $v_i$ is a \emph{descendent} of $v_j$ and similarly $v_j$ is an \emph{ancestor} of $v_i$. For some fixed graph $G$, natural number $l$, and $H \in IIM_l(G)$ let $H_i$ for $0 \le i \le l$ denote the induced subgraph of $G$ containing all vertices in levels 0 through $i$.
For the sake of brevity, we will often refer to graphs that can be generated by this model as \emph{IIM graphs}. 


In our study of this model, we found that many structural properties hold for all IIM graphs, regardless of the specifics of a randomized model. Thus, the focus of this paper is strictly the general deterministic model, to highlight the fundamental structural properties of complex network models using transitivity in their generation.
This allows us to focus on what is true about all iterated models, thus removing non-interesting questions from any future probabilistic analysis. The natural probabilistic model, which we discuss more in the conclusion, contains substantial edge dependencies, thus differentiating this model significantly from the standard Erd\H{o}s-R\`{e}nyi random graph. However, these dependencies indicate that specific analysis for individual applications will likely prove more fruitful than a general probabilistic analysis. By focusing on the general deterministic model, our results immediately extend all probabilistic cases, including those cases where a probabilistic analysis is unmanageable due to the intricacies of the edge dependencies.

\section{Spectral Gap}\label{spectral}

Complex networks will often exhibit large spectral gaps, and this is called Good Expansion \cite{estrada2006spectral}. This property creates sparse networks with patches of well-connected sub-graphs, and shows the clustered but non-uniform nature of complex networks \cite{estrada2006spectral, estrada2005complex}.
We begin this section by reviewing the definition of the spectral gap with respect to the eigenvalues of the normalized Laplacian. Then, we show that IIM graphs exhibit a spectral gap bounded away from zero, which is sufficient to demonstrate the property of Good Expansion. We did not attempt to optimize the lower bound and assuredly one could improve it slightly, but at the cost of messier computations.

Let $G$ be a graph on $n$ vertices; the adjacency matrix $A$ is the $n\times n$ matrix whose rows and columns are indexed by the vertices of $G$ and whose $ij^{th}$ entry is 1 if the corresponding vertices are adjacent and 0 otherwise. Let $D$ denote the diagonal degree matrix of $G$, that is the $n\times n$ diagonal matrix whose rows and columns are indexed by the vertices of $G$ and whose $ii$th entry is the degree of the corresponding vertex. For a graph $G$, with no isolated vertices, and adjacency matrix $A$ and diagonal degree matrix $D$, we define the normalized Laplacian matrix of $G$ as

\begin{equation*}
    L = I-D^{-1/2}AD^{-1/2}.
\end{equation*}

Note that the the eigenvalues of $L$ satisfy the following bounds $0 = \lambda_0 \le \lambda_1 \le \cdots \le \cdots \lambda_{n-1} \le 2$ (see e.g. \cite{spectralbook}). Thus the spectral gap of $L$ is defined by 
\begin{equation*}
    \lambda = \max\{|\lambda_1 -1|, |\lambda_{n-1}-1|\}.
\end{equation*}

We will use the following version of the Expander-Mixing lemma (see e.g. \cite{spectralbook}) to provide a lower bound on the spectral gap of any IIM graph with at least 4 levels.

\begin{lemma}\label{specmixlemma}
If $G$ is a graph with spectral gap $\lambda$, then, for all sets $X \subseteq V(G)$, 
\begin{equation}\label{specmixeq}
    \left| 2|E(X)| - \frac{(\vol(X))^2}{\vol(G)}\right| \le \lambda \frac{\vol(X)\vol(\overline{X})}{\vol(G)}.
\end{equation}
\end{lemma}

A more general version of the above lemma can be loosely interpreted as stating that graphs with nice, close to random, edge distribution have a small, close to zero, spectral gap. In particular, the spectral gap of the standard Erd\H{o}s-R\`{e}nyi random graph, above the connectivity threshold, is $o(1)$. More information can be found in \cite{purcilly_2020}. Thus by showing the spectral gap of IIM graphs are bounded away from 0 we confirm that the edge distribution of any IIM graph is not close to the standard uniform random models. This is more in line with many naturally occurring networks, such as social networks, further indicating that IIM graphs share more properties, and thus better mimic, social networks then traditional uniformly randomly generated graphs.

\begin{thm}
Let $G$ be any graph then for any $l \ge 4$ the spectral gap of any $H \in IIM_l(G)$ is at least $\frac{1}{15}$.
\end{thm}

\begin{proof}
Let $X$ be the set of vertices in level $l$ of $H$. As $X$ is an independent set $|E(X)|=0$, so (\ref{specmixeq}) simplifies to 
\begin{equation*}
    \frac{\vol(X)}{\vol(\overline{X})}\le \lambda.
\end{equation*}
Each vertex in $X$ is a copy of a vertex in $\overline{X}$; thus for $v \in X$ let $\overline{v} \in \overline{X}$ denote the precopy of $v$. Again as $X$ is an independent set we have $$\vol(\overline{X}) = 2|E(\overline{X})|+\vol(X).$$ 
We will show the middle inequality in
\begin{equation}\label{specfinal}
2|E(\overline{X})| = \sum_{\overline{x} \in \overline{X}}\deg(\overline{x}) \le 14 \sum_{x \in X} \deg(x) = 14 \vol(X)
\end{equation}
giving the desired bound $\frac{1}{15} \le \frac{\vol(X)}{\vol(\overline{X})}\le \lambda.$


As we will be comparing the degrees of vertices in $\overline{X}$ with their copies in $X$ recall that for each vertex $v \in X$:
\[ \deg(v)=
\begin{cases}
\deg_{\overline{X}}\overline{v}+1 & \textrm{if v is an clone of $\overline{v}$}\\
|X|-\deg_{\overline{X}}\overline{v}-1 & \textrm{if v is an anticlone of $\overline{v}$}.
\end{cases}\]
Let $Y \subset \overline{X}$ denote the set of vertices in $H$ in levels $l-1$ and $l-2$ and let $Z = \overline{X}\setminus Y$, and note that $|Y|= \frac{3|X|}{4}$ and $|Z|= \frac{|X|}{4}$. Further partition $Y$ and $Z$ into $Y = Y^c \cup Y^a$ and $Z = Z^c \cup Z^a$ to denote the vertices in each set that were cloned ($Y^c$ and $Z^c$) and anticloned ($Y^a$ and $Z^a$) in $X$. Furthermore, for any $U \subset \overline{X}$ let $A(U)$ be the set of vertices in $X$ which are anticlones of a vertex in $U$ and let $C(U)$ be the set of vertices in $X$ which are clones of a vertex in $U$.

\begin{figure}[h]
    \centering
    \includegraphics[height=3cm]{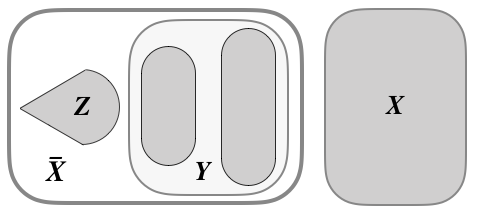}
    \caption{The subsets $X$, $\overline{X}$, $Y$, and $Z$ as defined within the IIM graph $H$. }
    \label{fig:specgapxyz}
\end{figure}

Recall that our goal is to balance the degrees within $\overline{X}$, of vertices in $\overline{X}$, with the total degrees of their copies in $X$ as in (\ref{specfinal}). It will be helpful to think of the vertices in $\overline{X}$ and $X$ in pairs $\overline{v}, v$ of vertices and their respective copies. As the degrees of $\overline{v}$ and $v$ are almost identical when $v$ is a clone, this case poses no difficulty. Furthermore, the degree of an anticlone is inversely proportional to that of its precopy. Note it only helps us if the anticlone has larger degree than its precopy; in fact, we can even easily afford anticlones with degrees that are only a seventh of their precopy given the inequality we're seeking. Thus we will only need to focus on pairs where the precopy of an anticlone has relatively large degree,  which for ease of computation will mean $\deg_{\overline{X}}(\overline{v}) \ge \frac{3|X|}{4}$.

No vertex in $Y$ can have relatively large degree, since $Y$ consists of an independent set of size $\frac{|X|}{4}$ and one of size $\frac{|X|}{2}$, thus bounding its largest degree to $\frac{3|X|}{4}$, which is small enough for our purposes. 
Therefore, our primary concern are vertices in $Z^a$ with large degree. We will handle this by dividing these ``problematic" vertices based on the number of neighbors they have in $Y^a$ versus $Y^c$. This allows us to pair the problematic vertices in $Z^a$ with ``well behaved'' vertices in $Y$ to show that as a whole the relevant degrees are balanced.

First, we consider the problematic vertices who have a sufficient number of neighbours in $Y^a$, that is vertices in $Z^a$ with large degree in $\overline{X}$ and many neighbours in $Y^a$. We define this set $Z^a_1$ as
\begin{equation*}
Z^a_1 = \left\{\overline{v} \in Z^a:  \deg_{\overline{X}}(\overline{v}) \ge \frac{3|X|}{4} \text{ and }\deg_{Y^a}(\overline{v}) \ge \frac{|X|}{4}\right\}.\end{equation*} 

If $Z^a_1$ is nonempty then $|Y^a| \ge \frac{|X|}{4}$. In this case we define $Y^a_1$ to be any subset of $Y^a$ of size $|Z^a_1|$, and note that $|Y^a_1|=|Z^a_1| \le |Z| = \frac{|X|}{4}$. As stated above every vertex $\overline{v} \in Y$ has $\deg_{\overline{X}}(\overline{v}) \le \frac{3|X|}{4}$ and thus their anticlones $v \in X$ must have $\deg(v) \ge \frac{|X|}{4}-1$.
We will also note that for these two bounds to be tight, $\overline{v}$ must be in level $l-2$ and a clone of an original vertex in $G$ among many other conditions. Thus we can simply use $|V(G)|$ to bound the number of vertices $v \in A(Y_1^a)$ with $\deg(v)=\frac{|X|}{4}-1$. Thus, using $|Z^a_1|=|Y^a_1|$ and $|\overline{X}|= |X|$ we have 
\begin{equation*}
    \sum_{\overline{v} \in Z^a_1} \deg_{\overline{X}}(\overline{v})+ \sum_{\overline{v} \in Y^a_1} \deg_{\overline{X}}(\overline{v}) \le |Z^a_1||X|+\frac{3|Y^a_1||X|}{4}=\frac{7|Z^a_1||X|}{4}
\end{equation*}
and 
\begin{equation*}
    \sum_{v\in A(Z^a_1)} \deg(v)+\sum_{v\in A(Y^a_1)} \deg(v) \ge \sum_{v\in A(Y^a_1)} \deg(v)\ge \frac{|Z^a_1||X|}{4} -|V(G)| \ge \frac{|Z^a_1||X|}{8}
\end{equation*}
where the final inequality uses the fact that $|X| \ge 8 |V(G)|$ since $l \ge 4$. We are potentially giving away a significant amount in each inequality; notice we are ignoring the contribution from vertices in $A(Z^a_1)$ entirely, but we seek to bound $\lambda$ away from zero rather than optimize. Therefore we have 
\begin{equation}\label{speceq1}
    \sum_{\overline{v} \in Z^a_1} \deg_{\overline{X}}(\overline{v})+ \sum_{\overline{v} \in Y^a_1} \deg_{\overline{X}}(\overline{v}) \le 14  \left( \sum_{v\in A(Z^a_1)} \deg(v)+\sum_{v\in A(Y^a_1)} \deg(v) \right)
\end{equation}

Now we consider the remaining problematic vertices, that have large degree in $\overline{x}$, but small degree in $Y^a$, which we define as follows 
\begin{equation*}
Z^a_2 = \left\{\overline{v} \in Z^a:  \deg_{\overline{X}}(\overline{v}) \ge \frac{3|X|}{4} \text{ and }\deg_{Y^a}(\overline{v}) < \frac{|X|}{4}\right\}. 
\end{equation*} 
Note that $\deg_{\overline{X}}(\overline{v}) \ge \frac{3|X|}{4}$ implies $\deg_Y(\overline{v}) \ge \frac{|X|}{2}$, as $|Z|= \frac{|X|}{4}$. Thus each $\overline{v} \in Z^a_2$ must have $\deg_{Y^c}(\overline{v}) \ge \frac{|X|}{4}$. Let $Y^c_2$ be the set of all vertices in $Y^c$ with neighbours in $Z^a_2$. We get the following bound on the size of these sets when $|Z^a_2| \ne 0$:
\[|Y^c_2| \ge \frac{|X|}{4} = |Z| \ge |Z^a_2|.\]
Therefore,
\begin{align}
\sum_{\overline{v} \in Z^a_2}\deg_{\overline{X}}(\overline{v}) + \sum_{\overline{v} \in Y^c_2} \deg_{\overline{X}}(\overline{v})&\le 4\sum_{\overline{v} \in Z^a_2} \deg_{Y^c_2}(\overline{v})+ \sum_{\overline{v} \in Y^c_2} \deg_{\overline{X}}(\overline{v}) \label{z2bound1}\\
& = 4\sum_{\overline{v} \in Y^c_2} \deg_{Z^a_2}(\overline{v})+ \sum_{\overline{v} \in Y^c_2} \deg_{\overline{X}}(\overline{v})\\
&\le  5 \sum_{\overline{v} \in Y^c_2} \deg_{\overline{X}}(\overline{v})\label{z2bound2}\\
&< 5\sum_{v \in C(Y^c_2)}\deg (v).\label{speceq2}
\end{align}
Here (\ref{z2bound1}) uses our assumption that each $\overline{v} \in Z^a_2$ must have $$\deg_{Y^c}(\overline{v}) \ge \frac{|X|}{4} \ge \frac{\deg_{\overline{X}}(\overline{v})}{4}.$$ 
Similar to the previous case, note that we do not need to consider the potential degrees of $A(Z^a_2)$ to obtain (\ref{speceq2}) as the degrees of $C(Y^c_2)$ alone are enough to provide a suitable bound. 

As mentioned earlier the remaining vertices in $Z^a \cup Y^a$ are easy to handle as they all have not too large degrees into $\overline{X}$. Specifically, every vertex $\overline{v} \in Z^a\setminus (Z^a_1 \cup Z^a_2)$ has $\deg_{\overline{X}}(\overline{v}) < \frac{3|X|}{4}$ so every vertex $v \in A(Z^a\setminus (Z^a_1 \cup Z^a_2))$ has $\deg(v ) \ge \frac{|X|}{4}$. Thus
\begin{equation}\label{speceq3}
    \sum_{\overline{v} \in Z^a \setminus (Z^a_1\cup Z^a_2)} \deg_{\overline{X}}(\overline{v})\;  \le\; \;  3 \sum_{v \in A(Z^a \setminus (Z^a_1\cup Z^a_2))} \deg(v)
\end{equation}
Similarly every vertex $\overline{v} \in Y^a$ has $\deg_{\overline{X}}(\overline{v}) < \frac{3|X|}{4}$ except potentially at most $|V(G)|$ vertices with $\deg_{\overline{X}}(\overline{v}) = \frac{3|X|}{4}$. Thus, again using the fact that $|V(G)| \le \frac{|X|}{8}$, we have
\begin{equation}\label{speceq4}
    \sum_{\overline{v} \in Y^a \setminus Y^a_1 } \deg_{\overline{X}}(\overline{v}) \le 8 \sum_{v \in A(Y^a \setminus Y^a_1)} \deg(v).
\end{equation}

\noindent The only vertices in $\overline{X}$ left to consider must have a clone in $X$ and thus 
\begin{equation}\label{speceq5}
    \sum_{\overline{v} \in \overline{X} \setminus (Z^a \cup Y^a \cup Y^c_2)} \deg_{\overline{X}}(\overline{v}) < \sum_{v \in C(\overline{X} \setminus (Z^a \cup Y^a \cup Y^c_2))} \deg(v).
\end{equation}
Combining equations (\ref{speceq1}), (\ref{speceq2}), (\ref{speceq3}), (\ref{speceq4}), (\ref{speceq5}) we have
\begin{equation*}
2|E(\overline{X})|=\sum_{\overline{v} \in \overline{X}} \deg_{\overline{X}}(\overline{v}) \le 14 \sum_{v \in X} \deg(v) = 14 \vol(X),
\end{equation*}
and thus we have the following bound on the spectral gap $$\frac{1}{15} \le\frac{\vol(X)}{2|E(\overline{X})| + \vol(X)} =  \frac{\vol(X)}{\vol(\overline{X})}\le \lambda,$$
as desired.

\end{proof}

In the previous result, we note that the only bound requiring the ratio of 14 is the one involving $Z^a_1$. Thus if $Z^a_1 = \emptyset$ our bound automatically improves from $\lambda \geq \frac{1}{15}$ to $\lambda \geq \frac{1}{9}$.

\section{Graph Properties}\label{graphprop}

The underlying graph structure within complex networks can yield interesting practical results \cite{easley2010networks, ilm, palu}. In \cite{fanbook}, four key properties of complex networks were identified, and continue to be a consensus of defining features: Large Scale, Evolution Overtime, Small World Properties, and Power Law Degree Distribution. Small world properties include low diameters and high clustering. In the following sections, we will determine the diameter of IIM graphs, which is in fact low. We then move on to discuss the domination number, the clique number, and end with a discussion on coloring and the chromatic number. The low domination of IIM graphs in conjunction with the non-zero spectral gap is an indicator of its well clustered behaviour. The domination number also can be used in information diffusion applications. Within our result on the clique number, and discussion of the chromatic number, we explore the overall optimization of fully independent and fully connected subsets of vertices. While gaining a better understanding of the general behavior of the clique and chromatic numbers is an interesting direction for future work, we also believe that given particular applications motivating additional restrictions on the IIM graphs of interest significantly more precise results can be obtained.

\subsection{Diameter}
Recall, that the \emph{diameter} of a graph $G$ denoted $\diam(G)$, is the maximum distance between any pair of vertices in $V(G)$, and the \emph{eccentricity} of a vertex $v \in V(G)$ denoted $\ecc_G(v)$, is the maximum distance from $v$ to any other vertex in the graph. As we will often be comparing distances in the initial graph $G$ with those in $H \in \mathcal{IIM}_l(G)$ we will use $\dist_G(u,v)$ for the distance, in $G$, between $u, v \in V(G)$ and simply $\dist(u,v)$ when referring to the distance in $H$.

In our following result, we give an upper bound for the diameter of a single iteration of the Iterated Independent Model. This result is extrapolated for all graphs in the subsequent corollary.

\begin{figure}[h]
    \centering
    \includegraphics[scale=.3]{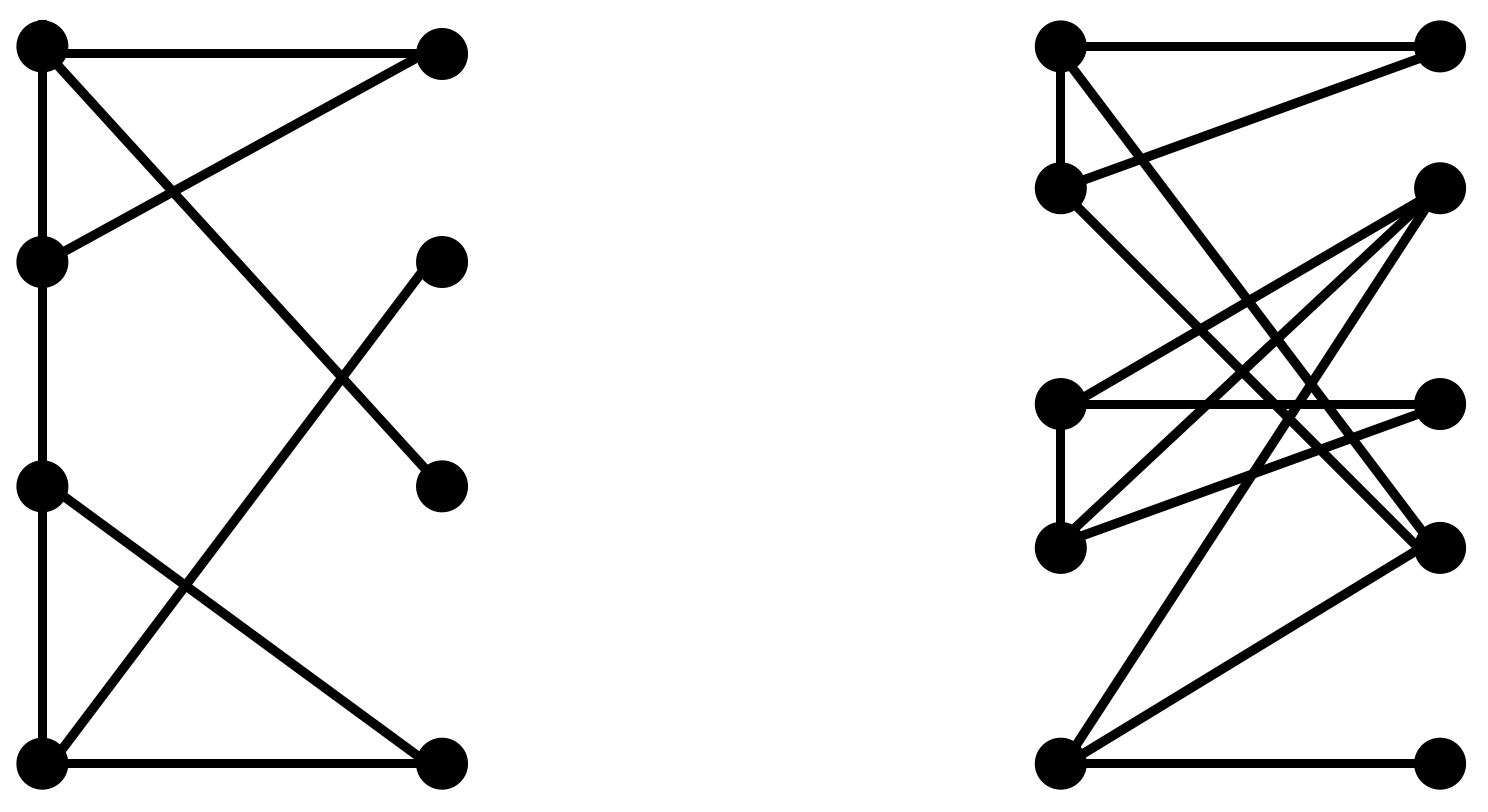}
    \caption{An illustration of the tight diameter bounds. Note that vertices in level 0 are on the left, with the copy of each vertex in level 1 at the same height on the right.
    }\label{diamex}
\end{figure}

\begin{thm}\label{diamthm}
Let $H$ be a connected graph in $IIM_1(G)$.
\begin{enumerate}
    \item If $G$ is connected then
    \begin{equation*}
        \diam(H) \le \max\{\diam(G), 5\}.
    \end{equation*}
    \item If $G$ is disconnected then
    \begin{equation*}
        \diam(H) \le 6.
    \end{equation*}
\end{enumerate}
\end{thm}
Before presenting the proof, we note that all the bounds are tight, and give an overview of a few observations. Figure~\ref{diamex} below on the left is an instance of a graph in $\mathcal{IIM}_1(P_4)$ with diameter $5$ and on the right is a graph in $\mathcal{IIM}_1(K_2 \cup K_2 \cup K_1)$ with diameter $6$. Additionally, note that if we clone every vertex the diameter of the resulting graph will always equal that of the starting graph.

\begin{proof}
We first prove part 1 of Theorem \ref{diamthm} and thus assume that $G$ is connected. We show that the diameter does not increase, except in the case where it increases to a maximum of 5.

Fix a connected $H \in \mathcal{IIM}_1(G)$. Assume that the diameter of $H$ has increased from that of $G$, since otherwise the result holds. 

\noindent \textbf{Claim 1.1: } $\diam(H) \le \diam(G)+2$. 

First we clearly have that if $v,w \in V(G)$ then $\dist_G(v,w) \le \dist(v,w)$. Consider next $v,w \in V(H)$ where exactly one of $v,w$ is in $V(G)$, say $v$ without loss of generality. 
When $w$ is a clone, $\dist(v,w)\le\dist_G(v,\overline{w})$, except of course when $v=w$ in which case $\dist(w,\overline{w})=1$. This is because $w$ is exactly adjacent to $\overline{w}$ and every neighbor of $\overline{w}$. Thus if $\overline{w},x, \ldots, v$ is a shortest path in $G$ between $\overline{w}$ and $v$ we have $w, x, \ldots, v$ is a path between $w$ and $v$ in $H$ of the same length. When $w$ is an anticlone, if $v,\overline{w}$ are not adjacent, then $\dist(v,w) = 1$; otherwise, consider some vertex $x \in V(G)$ not adjacent to $\overline{w}$, and thus adjacent to $w$, which must exist since $H$ is connected, then 
\[\dist(v,w) \le \dist(v,x) + 1 = \diam(G) +1.\]

Finally, when both vertices $v,w$ are in $V(H) \setminus V(G)$, they must each have at least one neighbour in $V(G)$, and thus the distance between $v,w$ is bounded above by $\diam(G)+2$. It is worth noting that together this implies that if the diameter increases from $G$ to $H$ then some new vertex $v\in V(H) \setminus V(G)$ must have $\ecc(v)=\diam(H)$.


\noindent\textbf{Claim 1.2:} For any $ u,v \in V(H)$ with at least one in $V(H) \setminus V(G)$, $\dist(u,v)$ can only be larger than $\diam(G)$ if $\dist(u,v) \le 5$. 

We separate the argument into cases based on the anticlone status of the vertices $ u,v$. \medskip

\noindent \underline{Case 1: Neither $u$ nor $v$ is an anticlone}

First assume $v$ is a clone of some vertex $\overline{v} \in V(G)$. Then for any $u \ne \overline{v} \in V(G)$ we have $\dist(u,v)\le \dist_G(u,\overline{v})$ and $\dist(v,\overline{v})=1$, as $v$ must be adjacent to every neighbor of $\overline{v}$.  Furthermore, if $u$ is a clone of $\overline{u}\in V(G)$ then, for the same reason as above $\dist(u,v)\le \dist_G(\overline{u}, \overline{v})$ unless $\dist_{G}(\overline{u}, \overline{v})=1$ in which case $\dist(u,v)=2$. Note that this now implies that if $\diam(H) > \diam(G)$ and $\diam(H)\ge 3$ we must have $\diam(H)=\ecc(v)$ for some $v$ an anticlone in $V(H) \setminus V(G)$. \medskip

\noindent\underline{Case 2: Exactly one of $u$ and $v$ is an anticlone}

Without loss of generality assume $v$ be such an anticlone and $\overline{v}$ be the precopy of $v$. First, for every $u \in V(G)$, consider distances between $v$ and $u$. If $\dist_{G}(u, \overline{v}) \ge 2$ then $u\overline{v} \notin E(G)$ so $uv \in E(H)$ and $\dist(u, v)=1$. If $\dist_{G}(u,\overline{v})= 1$ then let $z$ be any vertex in $G$ adjacent to $v$, which must exist since $H$ is connected and $V(H) \setminus V(G)$ is an independent set. We must then have that $\dist(u,v)\le \dist_{G}(u,z)+1$. However, if $\dist_{G}(u,z)\ge 4$ and $z,x, \ldots, u$ is some shortest path between $z$ and $u$ we must have $x\overline{v} \notin E(G)$, otherwise there would be a shorter path between $z$ and $u$ - namely $z,x, \overline{v},u$. This implies that $xv \in E(H)$ so $\dist(u,v) \le \dist_G(u,z)$. 

Similar to before if $u\in V(H)\setminus V(G)$ is a clone and $\overline{v}$ is non-adjacent to $\overline{u}$ then $\dist(u,v)= 2$. If this is not the case, i.e. $\overline{v}$ is adjacent to $\overline{u}$, then let $z, w_1, \ldots, w_k, \overline{u}$ be some shortest path between $z$ and $\overline{u}$, where we still assume $z$ is a neighbor of $v$ in $V(G)$. 
If $v$ is adjacent to any $w_i$ we have $\dist(u,v) \le \dist_G( \overline{u},z)$. If $v$ is not adjacent to any $w_i$, then $\overline{v}$ must be adjacent to all $w_i$, so a shortest path, between $z$ and $\overline{u}$ in $G$, must precisely be $z, w_1, \overline{v}, \overline{u}$, and $\dist(u,v) \le \dist_G(\overline{u},z) + 2 =5$. \medskip


\noindent \underline{Case 3: Both $u$ and $v$ are anticlones} 

As always let $\overline{u},\overline{v} \in V(G)$ be the precopies of $u$ and $v$, respectively. Assume $ua \in E(H)$ and $vb \in E(H)$ for $a, b \in V(G)$. As shown in Claim 1.1 $\dist(u,v) \le \dist_{G}(a,b)+2$. Assume $\dist(u,v) \ge 6$ and thus $\dist_{G}(a,b) \ge 4$. Let $a, w_1, w_2, \ldots, w_k, b$ be some shortest path between $a$ and $b$ in $G$. Note that if any of the following hold then $\dist(u,v) \le \dist_G(a,b)$:
\begin{enumerate}
    \item $uw_i \in E(H)$ for some $i \ge 2$ or $ub \in E(H)$.
    \item $vw_i \in E(H)$ for some $i \le k-1$ or $va \in E(H)$.
    \item $uw_1 \in E(H)$ and $vw_k \in E(H)$.
\end{enumerate}
If we assume, for that sake of contradiction, that $\dist(u,v) > \dist_G(a,b)$ then we may assume, without loss of generality, that $uw_i \notin E(H)$ for all $i$ and $ub \notin E(H)$ and thus since $u$ is an anticlone of $\overline{u}$ we have $\overline{u}w_i \in E(G)$ for all $i$ and $\overline{u}b \in E(G)$. However, this implies that $a, w_1, \overline{u}, b$ is a shorter path between $a$ and $b$ in $G$, a contradiction.  Thus it must be the case that either $\dist(u,v) \le 5$ or $\dist(u,v) \le \dist_G(a,b)$.

This concludes the first part of the theorem. For any possible pair of vertices in $V(H) \setminus V(G)$, we have that either $\dist(u,v) \le \diam(G)$ or that $\dist(u,v) \le 5$. Thus, whenever $G$ is connected, we have that $\diam(H) \le \max\{\diam(G), 5\}$.\bigskip

We now prove Part 2 of Theorem \ref{diamthm} where we assume $G$ is disconnected. First note that for $H \in \mathcal{IIM}_1(G)$ to be connected when $G$ is disconnected, there must be at least one anticlone in $V(H) \setminus V(G)$. We will split into two claims handling different possible configurations of $H$. 

\noindent \textbf{Claim 2.1:} If there is some precopy of an anticlone in $V(H) \setminus V(G)$ which does not dominate its component in $G$ then $\diam(H) \le 6$.

In this claim we assume there exists some anticlone $v \in V(H) \setminus V(G)$ such that if $\overline{v} \in V(G)$ is the precopy of $v$ then there is some $z \in V(G)$ such that $\dist_G(\overline{v},z)=2$. Let $y$ be a common neighbor of $\overline{v}$ and $z$. Note that this immediately implies $v$ is adjacent to at least one vertex in every component of $G$. We will now check the possible distances between any pair of vertices in $H$. \medskip

\noindent \underline{Case 1: $u, w \in V(G)$}

If $u, w \in V(G) \setminus N_G[\overline{v}]$ we have $\dist(u,w) \le 2$ as $u,v,w$ forms a path in $H$. Similarly if $u \in V(G) \setminus N_G[\overline{v}]$, but $w \in N_G[\overline{v}]$ we have $\dist(u,w) \le 5$ as $u,v,z,y,\overline{v},w$ forms a path in $H$. Finally if both $u, w \in N_G[\overline{v}]$ then they are within distance $2$ in $G$ and thus also in $H$. Thus for any $u, w \in V(G)$ we have $\dist(u, w) \le 5$. \medskip

\noindent \underline{Case 2: $u \in V(H) \setminus V(G)$ and $w \in V(G)$}

Since we are assuming $H$ is connected there is some $a \in V(G)$ adjacent to $u$. Thus we have $\dist(u,w) \le \dist(a,w)+1 \le 6$ as we have already shown $\dist(a,w) \le 5$ in Case 1. \medskip

\noindent \underline{Case 3: $u, w \in V(H) \setminus V(G)$ and both are clones}

As before let $\overline{u}\in V(G)$ be the precopy of $u$ and let $\overline{w} \in V(G)$ be the precopy of $w$. If $\overline{u}, \overline{w} \in N_G[\overline{v}]$ then $\dist(u,v) =2$ as $u,\overline{v},v$ is a path in $H$ (and $u$ and $w$ are in the same level and thus cannot be adjacent). If $\overline{u},\overline{w} \notin N_G[\overline{v}]$ then $\dist(u,w) \le 4$ as $u,\overline{u},v,\overline{w},w$ is a path in $H$. Finally, without loss of generality, assume $\overline{u} \in N_G[\overline{v}]$ but $\overline{w} \notin N_G[\overline{v}]$. Then $\dist(u,w) \le 6$ as $u,\overline{v},y,z,v,\overline{w},w$ is a path in $H$. \medskip

\noindent\underline{Case 4: $u, w \in V(H) \setminus V(G)$ and $u$ is a clone but $w$ is an anticlone}

Clearly if $\overline{u}\overline{w}\notin E(G)$ then $u,\overline{u},w$ is a path in $H$ so $\dist(u,w)=2$. Thus we may assume $\overline{u}\overline{w} \in E(G)$. Additionally, as we are assuming $H$ is connected there is some $c \in V(G)$ adjacent to $w$. If $\dist(\overline{u},c) \le 4$ then $\dist(u,w) \le 6$. If $\dist(\overline{u},c)>4$ then, by case 1, $\dist(\overline{u},c)=5$. Moreover, by the proof of Case 1 this implies either $c,v,z,y,\overline{v},\overline{u}$ or $\overline{u},v,z,y,\overline{v},c$ is a shortest path between $\overline{u}$ and $c$. If $c,v,z,y,\overline{v},\overline{u}$ is a path between $\overline{u}$ and $c$ then $\dist(u,w) \le 6$ as $w,c,v,z,y,\overline{v},u$ is a path of length $6$ in $H$. Now consider when $\overline{u},v,z,y,\overline{v},c$ is a path between $\overline{u}$ and $c$. If $w \overline{v}\in E(H)$ we would again have a path of length $6$, namely $u, \overline{u},v,z,y,\overline{v},w$, in $H$ between $u$ and $w$. If $w \overline{v}\notin E(H)$ then $\overline{w} \overline{v}\in E(G)$ as $w$ is an anticlone. Combining this with our earlier assumption that $\overline{u} \overline{w}\in E(G)$ we have $u, \overline{w}, \overline{v}, c, w$ is a path of length 4 between $u$ and $w$. \medskip

\noindent\underline{Case 5: $u, w \in V(H) \setminus V(G)$ and both are anticlones}

If $N_G[\overline{u}] \cup N_G[\overline{w}] \neq V(G)$ then $\dist(u,w)=2$. If this is not the case then $G$ must be comprised of only 2 components, one dominated by $\overline{u}$ and one dominated by $\overline{w}$. As $v$ is adjacent to at least one vertex in every component of $G$ let $a$ be a vertex in $\overline{u}$'s component adjacent to $v$ and let $b$ be a vertex in $\overline{w}$'s component adjacent to $v$. Then $u, b, v, a, w$ is a path of length $4$ between $u$ and $w$ in $H$. \bigskip

\noindent \textbf{Claim 2.2:} If the precopy of every anticlone in $V(H) \setminus V(G)$ dominates its component in $G$ then $\diam(H) \le 6$.

Our assumption means the neighborhood in $H$ of every anticlone is exactly the set of vertices in $G$ outside of it's precopies component (in $G$). In order for $H$ to be connected this means there must be at least three components in $G$ and at least two components containing a vertex anticloned in $H$. Let $u,v \in V(H) \setminus V(G)$ be such anticlones and let $\overline{u},\overline{v} \in V(G)$ be their precopies (i.e. $\overline{u}$ and $\overline{v}$ are in different components of $G$). Let $z$ be a vertex in $V(G)$ outside of $\overline{u}$ and $\overline{v}$'s components. 

We will now confirm that for any pair of vertices in $H$ their distance does not exceed 6. For any $a,b \in V(G)$ we must have that they are both non-adjacent to at least one of $\overline{u}$ and $\overline{v}$. Thus we know $\dist(a,b) \le 4$ as at least one of $a,u,b$, $a,v,b$, $a,u,z,v,b$, or $a,v,z,u,b$ is a path in $H$. If $a \in V(H) \setminus V(G)$ and $b \in V(G)$ then $\dist(a,b) \le 5$ since the connectedness of $H$ implies $ac\in E(G)$ for some $c \in V(G)$ and we just showed the distance between between $b$ and $c$ is at most $4$. Similarly if both $a, b \in V(H) \setminus V(G)$ then $\dist(a,b) \le 6$ as both $a$ and $b$ must have at least one neighbor in $G$.

\end{proof}
Theorem \ref{diamthm} immediately implies the following corollary. 
\begin{cor}
For any graph $G$, $l \in \N$, and connected $H \in \mathcal{IIM}_l(G)$ we have 
\begin{equation*}
    \diam(H) \le \max\{\diam(G), 6\}.
\end{equation*}
\end{cor}

\subsection{Domination number}

For a given graph $G$ we say $\mathcal{D} \subseteq V(G)$ is a dominating set of $G$ if for every $v \in V(G) \backslash \mathcal{D}$ there is a $u \in \mathcal{D}$ adjacent to $v$. We will use $\dom(G)$ for the size of the smallest dominating set(s) of $G$ and let $\Dom(G)$ denote such a minimum dominating set.  

Below we give two bounds on the dominating number of IIM graphs with different conditions on the initial graph. It is worth noting that if $H \in \mathcal{IIM}_1(G)$ is the graph formed by cloning every vertex then any dominating set of $G$ is a dominating set of $H$. We will expand this notion below, in what we call dual-dominating sets, to compensate for IIM graphs where some, but potentially not all, vertices are anticloned.
\begin{defn}
Given any graph $G$ if $\mathcal{D}$ is a dominating set of $G$ then $\mathcal{D}$ is a dominating set of every possible graph in $\mathcal{IIM}_1(G)$ if and only if every vertex in $G$ is also non-adjacent to at least one vertex in $D$. We will call such a set a \emph{dual dominating set} of $G$.
\end{defn}

\begin{thm}
Let $H \in \mathcal{IIM}_i(K_n)$ for some positive integers $n$ and $i$. Assume that the first anticlones in $H$ appeared at level $l$. Let $a_l(H)$ be the number of vertices in level $0$ anticloned in level $l$. 
\begin{equation*}
\dom(H) \le
    \begin{cases}
    4 & \text{ if }a_l(H)=0\\
    a_l(H)+3 & \text{ else}.
    \end{cases}
\end{equation*}

\end{thm}
\begin{proof}

For any graph $G$ and $H \in \mathcal{IIM}_1(G)$ let $\mathcal{D}$ be a dual-dominating set of $G$. By definition for every $v \in G$ there is some $u \in\mathcal{D}$ adjacent to $v$ and some $w \in \mathcal{D}$ non-adjacent to $v$. Thus the copy of $v$ in level 1 of $H$ will be adjacent to exactly one of $u$ and $w$. Therefore, $\mathcal{D}$ is also a dual dominating set of $H$, and moreover $\mathcal{D}$ will remain a dominating set for every graph in $\mathcal{IIM}_i(G)$ for every $i \ge 0$. 

Consider initializing with $K_1$. If in our first step we anticlone the initial vertex the two vertices, in level 0 and level 1, form a dual dominating set and thus any subsequent IIM graph will have domination number at most two. If we clone the initial vertex in the first step then this is equivalent to initializing with a $K_2$. Thus, moving forward, we may assume $H \in \mathcal{IIM}_i(K_n)$ for $n \ge 2$.

Let $v_0 \in V(K_n)$ and let $A = \{v_0, v_1, \ldots, v_{l-1}\}$ where each $v_j$ is the clone, in level $j$, of $v_{j-1}$. Notice that this exists as all levels up to level $l$ contain only clones and it forms a complete subgraph of $H$. 
Further, we designate some vertex in level 1 that is not $v_1$, to be a special vertex $x$. \bigskip


Recall, $a_l(H)$ is the number of vertices of the original graph $G$ that were cloned in level $l$. We continue the proof by separating into two cases, where either no original vertices are anticloned in level $l$, or at least one is anticloned. \medskip

\noindent \underline{Case 1: $a_l(H)=0$} 

Let $y$ be some anticlone appearing in level $l$. We claim that $$\mathcal{D}= \{ v_0, v_{l-1}, x, y\}$$ is a dual dominating set of $H_l$ and thus a dual dominating set of $H$. In fact $\{v_0, v_{l-1}, x\}$ is a dominating set of $H_l$, but we require the additional vertex $y$ to ensure the non-adjacency property. First note that $N_{H_{l-1}}(v_{l-1}) = A \cup V(H_0)$, and as no vertices in level 0 are anticloned by level $l$ we know $x$ is not only in level 1, but is also a clone of a level 0 vertex and thus $N_{H_{l-1}}(x)\cap A= v_0$. To confirm the non-adjacency property note that every vertex in $V(H_{l-1})$ is adjacent to every vertex in $V(H_0)$, and thus $y$ must be non-adjacent to every vertex in $V(H_0)$. Therefore we have,

\begin{itemize}
    \item $v_0$ dominates all of $V(H_{l-1})$ and all of the clones appearing in level $l$;
    \item $v_{l-1}$ is non-adjacent to $V(H_{l-1}) \setminus (A \cup V(H_0))$, and thus 
    dominates all anticlones of $V(H_{l-1})\setminus (A \cup V(H_0))$ in level $l$;
    \item $x$ is non-adjacent to $A \setminus \{v_0\}$, and thus dominates all anticlones of $A \setminus  \{v_0\}$;
    \item $y$ is non-adjacent to $V(H_0)$ and all vertices in level $l$.
\end{itemize}
\noindent Since we are assuming all vertices in $V(H_0)$ were cloned at level $l$ we have that $\mathcal{D}$ is a dual dominating set of $H_l$ and thus also of $H$, as desired. \medskip

\noindent \underline{Case 2: $a_l(H)\ge 1$.}

Note that if $u$ is an anticlone at level $l$ of some $v \in V(H_0)$ then $u$ is an isolated vertex in $H_l$. Let $U$ be the set of all such vertices (so $|U|=a_l(H)$). We claim the set $$\mathcal{D}=\{v_0, v_{l-1}, x\} \cup U$$ is a dual dominating set of $H_l$ and thus a dual dominating set of $H$. First we confirm that $\mathcal{D}$ is a dominating set of $H_l$. We have,

\begin{itemize}
    \item $v_0$ dominates every vertex in $V(H_{l-1})$ and all of the clones appearing in level $l$;
    \item $v_{l-1}$ dominates all anticlones of $V(H_{l-1})\setminus (A \cup V(H_0))$;
    \item $x$ dominates all anticlones of $A \setminus  \{v_0\}$;
    \item $U$ dominates itself (all anticlones of $V(H_0)$).
\end{itemize}
To confirm the non-adjacency property we simply need to note that since all $u \in U$ are isolates in $H_l$ we have $\{v_0, u\}$ is an anti-dominating set of $H_l$ for any $u \in U$. Thus $\mathcal{D}$ is a dual dominating set of $H_l$ and thus a dual dominating set of $H$, as desired. 
\end{proof}

In Theorem \ref{otherdom} we provide a bound on IIM graphs whose starting graph has dominating number at least 2. We leave open the question of a meaningful bound for IIM graphs who originate with a graph which is not a clique, but does have a single dominating vertex. 

\begin{thm}\label{otherdom}
Let $G$ be any graph with $\dom(G) \ge 2$ and define 
\begin{equation*}
b(G)=\min_{u,v \in V(G), uv \notin E(G)}|(N_{G}(v) \cap N_G(u))|.
\end{equation*}
If $H \in \mathcal{IIM}_i(G)$, for any $i \in \N$
then 
\begin{equation*}
\dom(H) \le \dom(G)+b(G)+3.
\end{equation*}
\end{thm}
\begin{proof}
Identify some $u,v \in V(G)$ satisfying the minimality of $b(G)$, that is $uv \notin E(G)$ and $|(N_{G}(v) \cap N_G(u))|=b(G)$. Again we will assume the first anticlones in $H$ appear in level $l$. Let $A= \{v_0, \ldots, v_{l-1}\}$ where $v_0 = v$ and $v_j$ is the clone of $v_{j-1}$ in level $j$, and let $y$ be some anticlone in level $l$. It may be the case that $y$ is the anticlone of some vertex in $A$, but is not necessary nor a restriction. 

Let $Z(G)$ be a minimal set of vertices in $G$ that contains at least one non-neighbour of each vertex in the aforementioned intersection, $N_{G}(v) \cap N_G(u)$; specifically, such that for every $w \in N_{G}(v) \cap N_G(u)$ there is some $z \in Z(G)$ such that $wz \notin E(G)$. We know such a set exists in $G$ as $\dom(G) \ge 2$. We claim that $\mathcal{D}=\Dom(G) \cup Z(G) \cup \{u, v_{l-1},y\}$ is a dual dominating set of $H_l$ and thus $H$. Note that since $\Dom(G)$ dominates $V(H_{l-1})$ there is some vertex $w \in \Dom(G)$ such that $w$ is adjacent to $y$'s precopy and thus non-adjacent to $y$. Therefore we have,

\begin{itemize}
    \item $\Dom(G)$ dominates every vertex in $V(H_{l-1})$ as well as all of the clones appearing in level $l$;
    \item $v_{l-1}$ is non-adjacent to $V(H_{l-1}) \setminus (A \cup N_G(v))$, and thus dominates all anticlones of $V(H_{l-1})\setminus (A \cup N_G(v))$;
    \item $u$ is non-adjacent to $A \cup (N_G(v) \setminus N_G(u))$, and thus dominates all anticlones of $A \cup (N_G(v) \setminus N_G(u))$;
    \item $Z(G)$ dominates all the anticlones of $N_{G}(v) \cap N_G(u)$;
    \item for each $w \in N_G(v) \cap N_G(u)$ there is a $z \in Z(G)$ such that $w$ is non-adjacent to $z$;
    \item $y$ is non-adjacent all vertices in level $l$;
    \item there is some $w \in \Dom(G)$ which is non-adjacent to $y$.
\end{itemize}

Thus $\mathcal{D}$ is a dual dominating set of $H_l$ and thus also $H$ of size at most $\dom(G)+b(G)+3$, as desired.

\end{proof}

\subsection{Clique number and Coloring}

Given a graph $G$, some integer $l \in \N$, and $v \in V(G)$ then for any $H \in \mathcal{IIM}_l(G)$ the set of vertices containing $v$ and all it's clones form a clique. Furthermore, all anticlones of $v$ also form a clique. Thus we trivially have a lower bound of  $\lceil (l+1)/2\rceil$ for the clique number. However, we significantly improve this lower bound in Theorem \ref{cliquebound} by more carefully examining what it takes to avoid increasing the size of a clique when adding a new level. We will use $\omega(G)$ to denote the clique number of $G$.
\begin{defn}\label{nonadj}
Given any graph $G$ we say we have a \emph{non-adjacent triple} if there are three disjoint subsets of vertices under some ordering $A_1, A_2, A_3 \subseteq V(G)$ such that there exists

\begin{enumerate}
    \item $v_1 \in A_1$ non-adjacent to every vertex in $A_2$;
    \item $v_2 \in A_2$ non-adjacent to every vertex in $A_3$;
    \item $v_3 \in A_3$ non-adjacent to every vertex in $A_1$.
\end{enumerate}
\end{defn}

\begin{figure}[h]
    \centering
    \includegraphics[scale=.85]{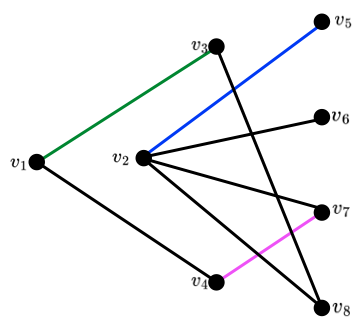}
    \caption{An example of a non-adjacent triple. Here $A_1$ is shown in green, $A_2$ is shown in pink, and $A_3$ is shown in blue.
    }\label{cliquefig}
\end{figure}

In Figure~\ref{cliquefig}, we give an example of a particular graph $H \in \mathcal{IIM}_3(K_1)$ with clique number 2. Here level 0 is $v_1$ and the vertices in level $i$, for $i\ge 1$ use labels $\{v_{2^{i-1}+1}, \ldots, v_{2^{i}}\}$ where vertex $v_{2^{i-1} +j}$ is the copy, in level $i$, of vertex $v_j$. Note that $A_1=\{v_1,v_3\}$, $A_2 = \{v_4,v_7\}$, and $A_3= \{v_2,v_5\}$ form a non-adjacent triple of $K_2$'s. In this instance $v_3$ is non-adjacent to both $v_4$ and $v_7$, $v_4$ is non-adjacent to both $v_2$ and $v_5$, and both $v_2$ and $v_5$ are non-adjacent to $v_1$ and $v_3$. 

We will show in the proof below that given any graph $G$ with a non-adjacent triple of $K_2$'s, any graph in $\mathcal{IIM}_3(G)$ must have clique number of at least 4. This is because, for any $H \in \mathcal{IIM}_3(G)$, at least one of the $A_i$'s will be contained in a $K_4$ regardless of the clone/anticlone choices. Moreover, one can ensure that the clique number is exactly $4$ in a particular graph in $\mathcal{IIM}_3(G)$ by carefully choosing which vertices are cloned/anticloned. This will be explained in more generality after the statement and proof of Theorem \ref{cliquebound}.

\begin{thm}\label{cliquebound}
For any graph $H\in \mathcal{IIM}_k(K_1)$ and $k \in \mathbb{N}$, 
\[\omega(H) \geq 2 + \left\lceil \frac{2(k-5)}{3} \right\rceil. \]
\end{thm}
\begin{proof}
Recall that $H_0 = K_1$, and the induced subgraph of $H$ containing vertices in levels $0$ through $i$ is denoted $H_i$. Furthermore, it must be true that $H_i \in \mathcal{IIM}_1(H_{i-1})$ since we form $H_i$ by creating one level of copies from $H_{i-1}$. 
Assume for some $l$, $H_l$ contains a non-adjacent triple $A_1, A_2, A_3 \subseteq V(H_l)$  such that $H_l[A_i]$ are cliques for all $i \in [3]$. In the remaining proof, we will use the notation $v_i$ to denote a specific vertex in $A_i$ witnessing the non-adjacency with $A_{i+1 (\mod 3)}$, as in Definition \ref{nonadj}. Note that as $v_i \in A_i$ and is non-adjacent to every vertex in $A_{i+1 (\mod 3)}$ we know that if $v'_i$ is the copy of $v_i$ in level $l+1$ then either $A_i \cup v'_i$ or $A_{i+1 (\mod 3)} \cup v'_i$ forms a clique in $H_{l+1}$. \bigskip 

\noindent{\bf Claim:} Any combination of cloning or anticloning the $v_i$'s results in a non-adjacent triple, $B_1, B_2, B_3$ in $H_{l+1}$ where 
\begin{enumerate}
    \item $A_i \subseteq B_i$ for all $i$;
    \item each $B_i$ induces a clique in $H_{l+1}$;
    \item for at least two distinct $i, j \in [3]$ we have $|B_i|=|A_i|+1$ and $|B_j|=|A_j|+1$. 
\end{enumerate}
\noindent \underline{Case 1:} At least two of the $v_i$ are cloned in level $l+1$.\\
Without loss of generality assume $v_1$ and $v_2$ are cloned and let $v'_1$ and $v'_2$ denote their clones in level $l+1$. Then let $B_1=A_1 \cup v'_1$, $B_2 = A_2 \cup v'_2$ and $B_3 = A_3$. Clearly we have the three required conditions and the $B_i$'s must be non-adjacent as for each $i$ we still have $v_i$ is non-adjacent to every vertex in $B_{i+1 \mod 3}$. \medskip

\noindent \underline{Case 2:} Exactly one $v_i$ is cloned.\\
Without loss of generality assume $v_1$ is cloned and $v_2$ and $v_3$ are anticloned. Let $v'_1$ denote the clone of $v_1$ in level $l+1$ and let $v''_2$ and $v''_3$ denote the anticlones of $v_2$ and $v_3$ in level $l+1$. Let $B_1 = A_1 \cup v'_1$, $B_2 = A_2$, and $B_3 = A_3 \cup v''_2$. In this case as $v''_2$ must be non-adjacent to every vertex in $A_2$ and adjacent to every vertex in $A_3$ we again clearly satisfy the three requirements and each $v_i$ is still non-adjacent to every vertex in $B_{i+1 (\mod 3)}$. \medskip

\noindent \underline{Case 3:} Every $v_i$ is anticloned.\\
Let $v''_i$ denote the anticlone of $v_i$ in level $l+1$. Let $B_1 = A_1 \cup v''_3$, $B_2 = A_2 \cup v''_1$, and $B_3 = A_3 \cup v''_2$. Given that each $v_i$ is adjacent to every vertex in $A_i$ and non-adjacent to every vertex in $A_{i+1 (\mod 3)}$ $B_3, B_2, B_1$ forms a non-adjacent triple witnessed by $v_2'', v_1'', v_3''$. Specifically, in this case we have $v''_3$ is in $B_1$ and non-adjacent to $B_3$, $v''_1$ is in $B_2$ and non-adjacent to $B_1$, and $v''_2$ in in $B_3$ and non-adjacent to $B_2$. \medskip 

This completes the proof of the claim. 
Thus if there exists a non-adjacent triple $A_1, A_2, A_3 \subseteq V(H_l)$ where each $A_i$ induces a clique and $|A_i| \ge n$ for all $i$. Then for all $j \ge 1$ $$\omega(H_{l+j}) \ge n + \left \lceil \frac{2j}{3}\right \rceil. $$ 

All that remains to verify is that in any $H \in \mathcal{IIM}_4(K_1)$ we have a non-adjacent triple of $K_2$'s. This can be seen by solely considering the edges between the four vertices, $v_1, v_2, v_3,$ and $v_4$ in level 3 and their copies in level 4. Note that if there are three vertices, w.l.o.g., $v_1, v_2, $ and $v_3$ in level 3 and that are cloned to form $v'_1, v'_2, $, and $v'_3$ in level 4 then we have $A_i = \{v_i, v'_i\}$ forms a non-adjacent triple with each subset inducing a $K_2$. If two vertices, w.l.o.g., $v_1$ and $v_2$ in the level 3 are cloned to form $v'_1$ and $v'_2$ in level 4 while one vertex, w.l.o.g, $v_3$ is anticloned to form $v''_3$ then $A_1 = \{v_1, v'_1\}$, $A_2= \{v_2, v'_2\}$, and $A_3= \{v_4, v''_3\}$ again forms a non-adjacent triple with each subset inducing a $K_2$. Finally if there are three vertices, w.l.o.g, $v_1, v_2, $ and $v_3$ in level 3 and that are anticloned to form $v''_1, v''_2, $, and $v''_3$ in level 4 then we have $A_1 = \{v_1, v''_2\}$, $A_2 = \{v_3, v''_1\}$, and $A_3 = \{v_2, v''_3\}$ forms a non-adjacent triple with each subset inducing a $K_2$. 

\end{proof}

The previous result gave an improved lower bound for the clique number of any IIM graph. However, we now explore the potential worst cases to provide better intuition as to how far the lower bound could be increased. While there are many cases where we can force the existence of a graph with a particular clique number, it is more challenging to bound the clique numbers of \emph{all} IIM graphs for any arbitrary $G$. The specific structures that occur in a particular proof of a bound, can interact across the course of growing large IIM graphs, so a tighter bound becomes more interesting with a particular application guiding the structure.

\begin{prop}
For sufficiently large $l \in \mathbb{N}$, there is some $H\in \mathcal{IIM}_l(K_1)$, such that $$\omega(H) \le l- \log^*(l).$$
\end{prop}
\begin{proof}


We say that sets $A_1, \ldots, A_k \subseteq H_l$ are \emph{a non-adjacent tuple} if for each $i \neq j$ there is a vertex $v_{i,j} \in A_i$ non-adjacent to every vertex in $A_j$. We say two sets, $U_1$ and $U_2$ are \emph{adjacent} if for all $u \in U$ we have a $v \in V$ such that $uv \in E(H)$ and for all $v \in V$ there is some $u \in U$ such that $vu \in E(H)$. Let $\{A_1, \ldots, A_k\}$ be a maximal\footnote{\label{note1}i.e. there is no other clique $A_{k+1}$ of size $\omega(H_l)$ that is pairwise non-adjacent to all the other $A_i$'s } non-adjacent tuple of cliques of size $\omega(H_l)$ in $H_l$. Then there is a sequence of clone/anticlone choices from $H_l$ forming a $H_{l+k} \in \mathcal{IIM}_k(H_l)$ where $\omega(H_{l+k}) =\omega(H_i)+k-1$.
This claim is equivalent to stating that if $\{A_1, \ldots, A_k\}$ is a maximal non-adjacent tuple of cliques of size $\omega(H_l)$ in $H_l$ then if $k \ge 2$ we have
\begin{enumerate}
    \item $\omega(H_{l+1})=\omega(H_l)+1$ and
   \item the smallest maximal set of non-adjacent tuples of cliques of size $\omega(H_{l+1})$ in $H_{l+1}$ has size at least $k-1$.
\end{enumerate}

We present a `worst case' construction for $H$. First partition the maximum cliques in $H_l$ into sets $B_1, \ldots, B_k$ such that for all $j$
\begin{enumerate}
    \item $A_j \in B_j$ and
    \item all of the cliques in $B_j$ are adjacent. 
\end{enumerate}
 Form $H_{l+1}$ by anticloning all vertices in $B_1$ and cloning all vertices non-adjacent to $B_1$. In this case we will have a maximal set of non-adjacent tuples of cliques, $A^{l+1}_2, \ldots, A^{l+1}_k$, of size $\omega(H_{l+1})=\omega(H_l)+1$ in $H_{l+1}$ where $ A^{l+1}_j= A_j \cup v''_{1,j}$ for all $2 \le j \le k$ (where $v''_{1,j}$ is the anticlone of $v_{1,j}$ in level $l+1$). Similarly let $B^{l+1}_j$ be the set containing all $B_j$ as well as their neighbors in level $l+1$.
 Note that every maximum clique in $H_{l+1}$ must contain a maximum clique in $H_l$ which was in some $B_i$ with $i \ge 2$. Iterate this process where in level $l+2$ we anticlone all vertices in $B^{l+1}_2$ and clone all vertices non-adjacent to $B^{l+1}_2$.
 By continuing this process until you reach $H_{l+k}$ we guarantee that each maximum clique in $H_l$ grew in all but one level between $H_l$ and $H_{l+k}$. Furthermore, any maximum clique in $H_{l+k}$ must contain a clique in $H_l$ of size at least $\omega(H_l)-1$. Thus even having the crude bound of $|V(H_l)|/(\omega(H_l)-1)$ on the size of the smallest maximal set of non-adjacent tuples of cliques of size $\omega(H_{l+k})$ in $H_{l+k}$ provides an eventual bound of $\omega(H_l) \le l-\log^*l$.

\end{proof}

While it may be interesting to attempt to more closely control the size of the smallest maximal set of non-adjacent tuples of cliques as this process continues we believe this would be of most use given a specific context. 

It is a natural question to ask if there is a fixed relationship between the clique number and the chromatic number of any IIM graph. Studying the chromatic number of IIM graphs poses an interesting challenge. Let $H \in \mathcal{IIM}_k(G)$ for some fixed $k$ and graph $G$. Since each level in $H$ is an independent set we trivially have $\chi(H_{l+1}) \le \chi(H_l)+1$ for any $l$. 

With the goal of investigating the chromatic number of $H$ one might hope that any coloring of $H_l$ with $\chi(H_l)$ colors properly extends to a coloring of $H_{l+1}$ with $\chi(H_{l+1})$ colors. Here we make the distinction between extending a coloring, and giving a proper $\chi(H_{l+1})$ coloring, in that an extension of a coloring leaves all colors previously assigned to $V(H_l)$ intact, and chooses new colors only for the vertices in level $l+1$. In general a proper coloring would allow for a complete recoloring of all levels to minimize the colors. Unfortunately, in Example~\ref{colorex}, we show that restricting to only extensions is not necessarily optimal. In the subsequent analysis, we focus on the extension of colorings in order to explore the requirements based on the new vertices created.

\begin{defn}
Fix a graph $H$ and a coloring of the vertices with $c$ colors. 
\begin{enumerate}
    \item We say $v \in V(H)$ has a \emph{rainbow neighborhood} if for each color there is at least one vertex in $N[v]$ assigned that color. Similarly we say $v$ has a \emph{rainbow anti-neighborhood} if for every color there is at least one vertex in $V(H) \setminus N[v]$ assigned that color. 
    \item We say $H$ has a rainbow pair if least two vertices in $H$ have both a rainbow neighborhood and a rainbow anti-neighborhood.
\end{enumerate}

\end{defn}

\begin{prop}\label{color1}
Let $H \in \mathcal{IIM}_k(G)$ for some fixed $k$ and graph $G$. For some $l<k$ fix a proper coloring of $H_l$ with $c$ colors such that under this coloring $H$ has a rainbow pair. Then any proper coloring of $V(H_{l+1})$ extending the fixed coloring of $H_l$ must use at least $c+1$ colors, and if such a coloring using exactly $c+1$ colors then it must also have a rainbow pair. 
\end{prop}

\begin{proof}
Let $v, w \in V(H_l)$ be a rainbow pair under the fixed coloring of $H_l$ with $c$ colors. Thus in $H_{l+1}$ both of their copies in level $l+1$ must be assigned the new color. Note that we only need to add a single new color since level $l+1$ forms an independent set. Furthermore as their copies are in the same level, and thus non-adjacent, we know that they now form a rainbow pair in $V(H_{l+1})$ regardless of how the rest of the level $l+1$ vertices are colored. 
\end{proof}

If it were true that any coloring of $V(H_l)$ with $\chi(H_l)$ colors extends to a coloring of $V(H_{l+1})$ with $\chi(H_{l+1})$ then if a coloring of $V(H_l)$ with $\chi(H_l)$ has a rainbow pair we would know $\chi(H_i) = \chi(H_l)+i-l$ for any $l \le i$. However, Example \ref{colorex} shows this is not the case.

\begin{example}\label{colorex}
Consider $\overline{K_2}$ which clearly has a chromatic number of 1. However, if we color $\overline{K_2}$ with only one color then the two vertices make up a rainbow pair. Thus any subsequent coloring of $H \in \mathcal{IIM}_k(\overline{K_2})$ where, for any $1 \le l \le k$, the coloring of $H_l$ extends the coloring of $H_{l-1}$ using the minimum possible number of colors for an extension, must use $l+1$ colors. If this was always best possible then for any $k$ and $H \in \mathcal{IIM}_k(\overline{K_2})$ we would have $\chi(H)=k+1$. However, Figure~\ref{coloringfig} demonstrates this is not the case.
\begin{figure}[h]
    \centering
    \includegraphics[scale=.45]{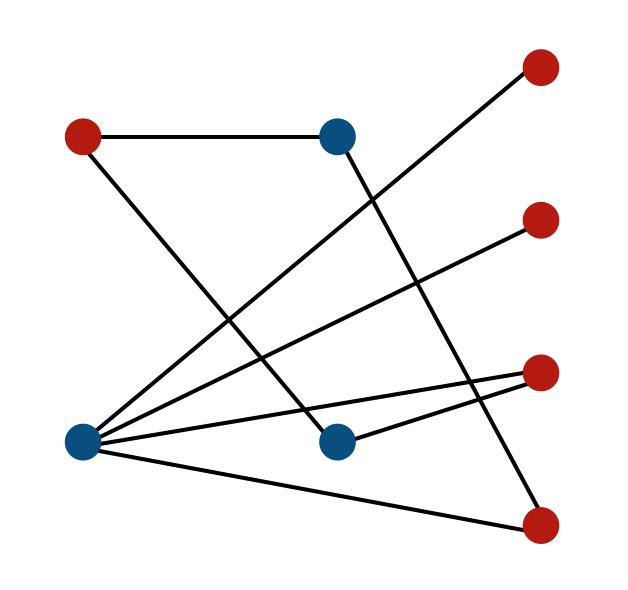}
    \caption{A specific $H \in \mathcal{IIM}_2(\overline{K_2})$ with a proper two coloring}
    \label{coloringfig}
\end{figure}

\end{example}

\section{Induced Subgraphs}\label{induced}

 In \cite{ilm} the authors showed the that for any fixed graph $F$ and any starting graph $G$ there is some, sufficiently large $k$, such that $F$ is an induced subgraph of every graph in $\mathcal{ILM}_k(G)$. While the proof in \cite{ilm} did rely on the full strength of the ILM class of graphs, namely knowing that for each level the vertices are either all clones or all anticlones, we are able to extend the result showing that for some $k$, potentially much larger than the $k$ needed in the ILM case, $F$ is an induced subgraph of every graph in $\mathcal{IIM}_k(G)$. Our result comes in two stages. First, in Lemma \ref{inducedlemma}, we give the result for IIM graphs with a very specific structure. Then, in Theorem \ref{inducedthm}, we show a method of constructing an induced copy of $F$ which, if it fails, immediately implies the specific structure described in Lemma \ref{inducedlemma}.

Given any graph $G$ and $v \in V(G)$ and $H \in \mathcal{IIM}_l(G)$, we begin by defining the following two subsets $S_c(v), S_a(v) \subset V(H)$:

 \begin{align*}
   S_c(v) = &\{u \in V(H): u=v \text{ or } u \text{ is a clone of }v\};\\
    S_a(v) = & \{u \in V(H):u \text{ is an anticlone of }v\}.
\end{align*}
Both of these sets induce complete subgraphs in $H$, and at least one contains at least $\lceil l/2 \rceil$ vertices. Thus for any $n$ there is a sufficiently large $l$ such that $\mathcal{IIM}_l(G)$ contains $K_n$ as an induced subgraph. Therefore, for the rest of the section we will restrict to simply showing that given any graph $F$ there is a induced copy of $F$ in $\mathcal{IIM}_k(K_n)$ for sufficiently large $k$ and $n = |V(F)|$. We begin with a useful special case, only considering clones, in Lemma \ref{inducedlemma}. In the full statement in Theorem \ref{inducedthm}, we are able to utilize this lemma or exploit the non-existence of the conditions.

\begin{lemma}\label{inducedlemma}
Let $F$ be any fixed graph on $n$ vertices and $m$ edges. For sufficiently large $k$ consider any $H \in \mathcal{IIM}_k(K_n)$. Assume there exists some $\mathcal{L} = \left\{l_0, \ldots, l_{{n \choose 2} - m}\right \}$ a subset of the indices of the levels,  and a corresponding set $U \subset V(H)$ such that 
\begin{enumerate}
    \item $l_0 =0$, and $l_{i+1} > l_i$ for all $i$;
  \item $U$ only contains vertices in levels indexed by $\mathcal{L}$;
  \item the entire vertex set of the initial $K_n$, i.e. level 0, appears in $U$;  
  \item for each $i \ge 1$ 
    \begin{enumerate}
       \item all of the copies in level $l_{i+1}$ of vertices in $U$ intersect level $l_j$ for $0\le j \le i$ are clones, and
      \item the set of vertices in $U$ intersect level $l_{i+1}$ is exactly the set of clones of the vertices lying in $U$ intersect level $l_j$ for $0\le j \le i$. 
    \end{enumerate}
\end{enumerate}
Then $G$ contains an induced copy of $F$.
\end{lemma}

To visualize these conditions think of building $U$ level by level. We start with all of level 0 then we assume that in level $l_1$ all the level 0 vertices are cloned and place all of those clones into $U$. Next we move to level $l_2$ where we assume all of the vertices thus far placed in $U$, i.e. level 0 and their clones in level $l_1$, are again cloned and again we place these clones in $U$. Continue this process until you reach level $l_{{n \choose 2}-m}$.

\begin{proof}
As before we use the notation $H_l$ for the induced subgraph of $H$ containing exactly the vertices in all levels up to, and including, level $l$. Label the vertex sets of $F$ and $K_n$ by $\{x_1, \ldots, x_n\}$. Let $ E(K_n) \setminus E(F) = \left\{e_1, \ldots, e_{{n \choose 2}-m}\right\}$.
We will progressively find induced copies of $K_n - \{e_1, \ldots, e_i\}$ in $H_{l_{i}}$.

For any $x_a$ in level 0 let $X_a = \left\{x^0_a, x^1_a, \ldots x^{{n \choose 2} - m}_a\right\}$ where $x^0_a=x_a$ and for any $i \ge 0$ we have $x^{i+1}_a$ is the vertex in level $l_{i+1}$ which is a clone of $x^{i}_a$. First note that if $e_1 = (x_a, x_b)$ then the subgraph of $H_{l_1}$ induced by $\left(\{x_1, \ldots, x_n\}\setminus \{x_a, x_b\}\right) \cup \{x^1_a, x^1_b\}$ is an induced copy of $K_n - e_1$. 

Now, for the sake of induction, assume we have an induced copy of $K_n - \{e_1, \ldots, e_i\}$ in $G_{l_{i}}$ containing exactly one vertex from each $X_a$; we will refer to this subgraph as $F_i$. Let $e_{i+1}=(x_c, x_d)$ and let $\{w\}= V(F_i) \cap X_c$ and $\{z\}= V(F_i) \cap X_d$. By the assumptions of our lemma both $w$ and $z$ have a clone, $w^{i+1}$ and $z^{i+1}$, in level $l_{i+1}$. Thus $(V(F_i) \setminus \{w,z\}) \cup \{w^{i+1}, z^{i+1}\}$ is an induced copy of $K_n - \{e_1, \ldots, e_{i+1}\}$ in $G_{l_{i+1}}$, as desired.

\end{proof}

Clearly not every graph in $\mathcal{IIM}_k(K_n)$ contains a $\mathcal{L}$ and $U$ satisfying the conditions of Lemma \ref{inducedlemma}. However, we show in Theorem \ref{inducedthm} that we can either build an induced copy of $F$ in a particular manner or such a $\mathcal{L}$ and $U$ must exist.

\begin{thm}\label{inducedthm}
Given any fixed graphs $F$ and $G$ there is a sufficiently large integer $k$ such that any $H \in \mathcal{IIM}_k(G)$ contains $F$ as an induced subgraph.
\end{thm}

\begin{proof}
Let $F$ be a fixed graph on $n$ vertices with $m$ edges. As we noted before Lemma \ref{inducedlemma} we only need to show that $F$ is an induced subgraph of any $H \in \mathcal{IIM}_k(K_n)$, for sufficiently large $k$. As before we label the vertex sets of $F$ and $K_n$ by $\{x_1, \ldots, x_n\}$. Similar to the proof of Lemma \ref{inducedlemma} we will show that we can find an induced subgraph by sequentially ``removing'' edges in $K_n$, until the desired $F$ is found. The removal is facilitated by replacing some vertices $x_i,x_j$ with copies $y_i, y_j$ such that $y_i$ and $y_j$ lie in the same level and satisfy the necessary adjacencies in $F$. We will show that if such $y_i$ and $y_j$ can never be found then we must have a $\mathcal{L}$ and $U$ satisfying the conditions of Lemma \ref{inducedlemma}.

For any $e = (x_i, x_j)$ we can find an induced copy of $K_n - e$ in $H_l$ if at level $l$ we see vertices, $y_i$ and $y_j$, which are descendants of $x_i$ and $x_j$, after an even number of anticloning steps. Specifically the graph induced by $(\{x_1, \ldots x_n\} \setminus \{x_i,x_j\})\cup \{y_i, y_j\}$ is an induced copy of $K_n \setminus e$ as $y_i$ and $y_j$ both appear in level $l$, and are thus non-adjacent, and their neighborhoods must contain $\{x_1, \ldots, x_n\}$ as an even number of anticloning steps were applied in each case. 

Now assume there is no such level. Then in every level either every descendent of $x_i$ or every descendent of $x_j$ must have been formed by an odd number of anticloning steps. Let $\mathcal{L}^q$ be the set of indices of the levels where every vertex in a such a level, who is a descendent of $x_q$ was formed by an odd number of anticloning steps. Additionally let $0 \in \mathcal{L}^q$ for every $q$. Clearly given a sufficiently large number of levels one of $|\mathcal{L}^i|$ and $|\mathcal{L}^j|$ is at least $n+{n \choose 2}-e(F)$.

Without loss of generality assume $|\mathcal{L}^i| \ge n+{n \choose 2}-e(F)$ and let $z$ be the anticlone of $x_i$ in the earliest nonzero index appearing in $\mathcal{L}^i$. Consider the set of vertices, $\mathcal{F}$, containing all vertices $v$ satisfying
\begin{enumerate}
    \item $v$ is a descendent of $z$, and
    \item any vertex $w$ which is both a descendent of $z$ and an ancestor of $v$ must be in a level whose index is in $\mathcal{L}^i$.
\end{enumerate}
As $z$ is an anticlone of $x_i$ every vertex in $\mathcal{F}$ must be formed only by cloning steps after $z$. Therefore let $f_1, \ldots, f_n$ be the set of vertices such that $z = f_1$ and each $f_i$ is the clone of $z$ appearing in the first level in $\mathcal{L}^i$ after the level $f_{i-1}$ appeared in. Again we have $\{f_1, \ldots, f_n\}$ induces a copy of $K_n$. Thus the induced subgraph of $H$ containing $\{f_1, \ldots, f_n\}$ and all of their clone descendants in a level indexed by $\mathcal{L}^i$ satisfies the conditions of Lemma \ref{inducedlemma} giving an induced copy of $F$ as desired.

Note that while we only showed how to remove a single edge, this process can be exactly repeated. Assume you have an induced copy of $F_r = F \setminus \{e_1, \ldots, e_r\}$. Proceed by only considering the subgraph of $H$ which is in $\mathcal{IIM}(F_r)$.
Specifically, if after removing $r$ edges your current vertex set, $\{a_1, \ldots, a_n\}$, contains only vertices in levels $l_0, l_{i_1}, \ldots, l_{i_r}$ then moving forward you only consider vertices which are descendants of an $a_i$ who both appear after level $l_{i_r}$ and, moreover, do not have any ancestors, except $a_i$, appearing in level $l_{i_r}$ or earlier.

\end{proof}

\section{Conclusion and future directions}

Overall, the explorations in this paper reveal underlying structures in graphs created through the IIM process. Thus, for any particular application, one can determine the desired properties, and restrict the IIM process in such a way to yield a theoretically accurate model. The insights created through specific applications of this model can be used to develop understanding of social network behaviours and information diffusion networks \cite{small2013information}.

Specifically, we introduced the IIM for social networks and were able to provide a vast generalization of many properties seen in the ILM, ILT, and ILAT models previously studied in \cite{ilm, ilat, ilt}. We proved the spectral gap of such graphs are guaranteed to be bounded away from zero, as seen in social networks; demonstrating that IIM graphs differ from the uniformly random graph model in edge density. In addition, we explored several graph properties showing bounds on the diameter, domination number, and clique number. We also began exploring coloring IIM graphs, demonstrating bounds when we restrict to only extending optimal colorings at every stage while also showing this is not guaranteed to produce an optimal coloring of the final graph.

There are still of course many open questions. It would be interesting to further explore coloring on these graphs specifically to gain an understanding of how far the chromatic number can be from the number of colors used when restricting to only extending the previous coloring at each level. Further investigation into the relationship between the chromatic number and clique number of these graphs is also of interest. Additionally, in \cite{ilm} the authors proved that the ILM becomes Hamiltonian after two non-consecutive anti-transitive steps. However, this questions becomes much more challenging in the IIM case, and it would be interesting to see what progress can be made towards classifying when IIM graphs are Hamiltonian.

For instance, given some graph $G$, natural number $l$, and $H \in \mathcal{IIM}_l(G)$ let $C$ be the set of vertices in $H_{l-1}$ cloned in level $l$ and let $A$ be the set of vertices in $H_{l-1}$ anticloned in level $l$. If, for some $k$, we can partition $C$ into $C_1, \ldots, C_k$ and $A$ into $A_1, \ldots, A_k$ such that
\begin{enumerate}
    \item each $C_i$ induces a Hamiltonian subgraph in $H_{l-1}$,
    \item each $A_i$ is an independent set with at least three vertices, and
    \item for each $i$ there exists some distinct $v,w \in C_i$, $u \in A_{i-1}$, and $x \in A_{i}$ such that $v u\notin E(H_{l-1})$ and $wx\in E(H_{l-1})$
\end{enumerate}
then $H$ must contain a Hamiltonian cycle. This is because if $C^c_i$ is the set of clones of $C_i$ in level $l$ then there is a Hamiltonian path in $C_i \cup C^c_i$ beginning with $v$ and ending with $w$ following the Hamiltonian path in $C_i$, but alternating with the clones. Specifically, if $v_1 v_2$ is an edge in a Hamiltonian cycle in $C_i$ and $v^c_1$ is the clone of $v_1$ in level $l$ then $v_1v^c_1$ and $v^c_1 v_2$ are consecutive edges in a Hamiltonian path of $C_i \cup C^c_i$. Similarly, as each $A_i$ is independent and has size at least three we can form a Hamiltonian path within each $A_i$ and their anticlones. Finally the third condition in the above list allows us to connect these Hamiltonian paths to form a Hamiltonian cycle in all of $H$. Unfortunately, it is not clear how to guarantee in general that such a partition must exits. Perhaps randomly forming $H$ as described below could provide some insight on when such a partition, and thus Hamiltonian cycle, does exist.

Finally, as mentioned in Section \ref{model}, we are interested in what can be said about a probabilistic version of this model. While this model can be randomized in many ways it is likely best to start with the following model. 
\begin{enumerate}
    \item Begin with a single node as $G_0 = K_1$.
    \item Define $G_t$ from $G_{t-1}$ as follows:
    \begin{itemize}
        \item for each node in $G_{t-1}$, clone the node with probability $\frac{1}{2}$ and anticlone the node with probability $\frac{1}{2}$.
    \end{itemize}
\end{enumerate}

Clearly we could further generalize this to the model where we clone with probability $p$ and anticlone with probability $1-p$ or ultimately allow the cloning probability to vary between vertices, perhaps based on the degrees. As this model randomly selects a member $\mathcal{H} \subseteq \mathcal{IIM}_l(K_1)$, any properties true above, will still hold. When particular structures are required in the presented results, determining the probability such a structure exists can help in determining the probability the property holds in general. An interesting component of this model is the significant edge dependencies. For instance, it is not too hard to show that any pair of edges $xy$ and $zy$ where $x$ and $z$ both appear in the same level as the pre-copy of $y$ are dependent. In general, however, even a complete analysis of the edge dependencies is not immediate. It would be interesting to first obtain a full classification of when two edges are dependent, and then potentially work towards a description of the degrees in such a model.

\pagebreak

\bibliographystyle{abbrv}
\bibliography{IIMbib}

\end{document}